\documentclass[letterpaper,12pt]{article}
\usepackage[margin=1in]{geometry}
\usepackage{amsmath}
\usepackage{amssymb}
\usepackage{amsthm}

\def\epsilon{\varepsilon}
\def\bigoh{\mathcal{O}}
\def\phi{\varphi}

\newtheorem{theorem}{Theorem}[section]
\newtheorem{lemma}[theorem]{Lemma}

\newtheorem{observe}[theorem]{Observation}
\newtheorem{remark1}[theorem]{Remark}

\newenvironment{remark}{\begin{remark1} \rm}{\end{remark1}}

\title{A fast algorithm for computing minimal-norm solutions
       to underdetermined systems of linear equations}
\author{Mark Tygert}
\date{September 8, 2009}

\begin{document}

\maketitle

\begin{abstract}
We introduce a randomized algorithm
for computing the minimal-norm solution
to an underdetermined system of linear equations.
Given an arbitrary full-rank matrix $A_{m \times n}$ with $m < n$,
any vector $b_{m \times 1}$, and any positive real number $\epsilon$
less than 1, the procedure computes a vector $x_{n \times 1}$
approximating to relative precision $\epsilon$ or better
the vector $p_{n \times 1}$ of minimal Euclidean norm satisfying
$A_{m \times n} \, p_{n \times 1} = b_{m \times 1}$.
The algorithm typically requires
$\bigoh( mn \, \log(\sqrt{n}/\epsilon) + m^3 )$
floating-point operations,
generally less than the $\bigoh(m^2 \, n)$ required
by the classical schemes based on $QR$-decompositions or bidiagonalization.
We present several numerical examples illustrating the performance
of the algorithm.
\end{abstract}

\section{Introduction}

Underdetermined systems of linear equations have arisen frequently
in modern statistics and data analysis, and have been attracting
much attention recently in various application domains;
see, for example,~\cite{bruckstein-donoho-elad}, \cite{candes-tao},
and~\cite{donoho}.
The solutions to underdetermined systems are not unique;
the present article focuses on the solutions whose Euclidean norms
are minimal.

Given a full-rank matrix $A_{m \times n}$ and a vector $b_{m \times 1}$,
with $m < n$, we would like to compute an accurate approximation
to the vector $p_{n \times 1}$ of minimal Euclidean norm satisfying
\begin{equation}
A_{m \times n} \, p_{n \times 1} = b_{m \times 1}.
\end{equation}
Classical algorithms using $QR$-decompositions or bidiagonalization require
\begin{equation}
\label{classical}
C_{\rm classical} = \bigoh(m^2 \, n)
\end{equation}
floating-point operations in order to compute $p_{n \times 1}$
(see, for example, Chapter~5 in~\cite{golub-van_loan}).

The present paper introduces a randomized algorithm that,
given any positive real number $\epsilon$ less than 1,
computes a vector $x_{n \times 1}$ approximating $p_{n \times 1}$
to relative precision $\epsilon$ or better
with respect to the Euclidean norm, that is,
the algorithm produces a vector $x_{n \times 1}$ such that
\begin{equation}
\| x_{n \times 1} - p_{n \times 1} \| \le \epsilon \, \| p_{n \times 1} \|,
\end{equation}
where $\| \cdot \|$ denotes the Euclidean norm.
This algorithm typically requires
\begin{equation}
\label{randomized}
C_{\rm rand} = \bigoh(mn \, \log(\sqrt{n}/\epsilon) + m^3)
\end{equation}
floating-point operations.
When $m$ is sufficiently large and $n$ is much greater than $m$
(that is, the system of linear equations is highly underdetermined),
then the cost in~(\ref{randomized}) is less than
the cost in~(\ref{classical}).
Moreover, in the numerical experiments of Section~\ref{numerical},
the algorithm of the present article runs substantially faster than
the standard methods based on $QR$-decompositions.

The present paper describes an algorithm optimized for the case
when the entries of $A_{m \times n}$ and $b_{m \times 1}$ are complex valued.
Needless to say, real-valued versions of our scheme are similar.
The present article has the following structure:
Section~\ref{notation} sets the notation.
Section~\ref{prelims} discusses a randomized linear transformation
which can be applied rapidly to arbitrary vectors.
Section~\ref{algorithm} describes the algorithm of the present paper.
Section~\ref{proofs} proves that the procedure succeeds with high probability.
Section~\ref{costs} estimates the computational costs of the algorithm.
Section~\ref{numerical} illustrates the performance of the scheme
via several numerical examples.
Section~\ref{conclusions} contains several concluding comments.

\section{Notation}
\label{notation}

In this section, we set notational conventions employed
throughout the present paper.

We abbreviate ``independent and identically distributed'' to ``i.i.d.''
We consider the entries of all vectors and matrices in this paper
to be complex valued.
For any vector $x$, we define $\|x\|$ to be the Euclidean ($l^2$) norm of $x$.
For any matrix $A$, we define $A^*$ to be the adjoint of $A$.
We define the condition number of $A$ to be the $l^2$ condition number of $A$,
that is, the greatest singular value of $A$ divided
by the least singular value of $A$.

For any positive integer $n$,
we define the discrete Fourier transform $F_{n \times n}$ to be the matrix
with the entries
\begin{equation}
\label{DFT}
F_{j,k} = \frac{1}{\sqrt{n}} \, e^{-2 \pi i (j-1) (k-1)/n}
\end{equation}
for $j,k = 1$,~$2$, \dots, $n-1$,~$n$, where $i = \sqrt{-1}$ and $e = \exp(1)$.

\section{Preliminaries}
\label{prelims}

In this section, we discuss a subsampled randomized Fourier transform.
\cite{ailon-chazelle}, \cite{drineas-mahoney-muthukrishnan-sarlos},
\cite{sarlos3}, and~\cite{sarlos4} introduced a similar transform
for similar purposes (these articles motivated us to write the present paper).

For any positive integers $l$ and $n$ with $l \le n$,
we define the $l \times n$ SRFT to be the random matrix
\begin{equation}
\label{SRFT}
T_{l \times n} = G_{l \times n} \, H_{n \times n},
\end{equation}
where $G_{l \times n}$ and $H_{n \times n}$ are defined as follows.

In~(\ref{SRFT}), $G_{l \times n}$ is the random matrix given by the formula
\begin{equation}
\label{SRFT2}
G_{l \times n} = S_{l \times n} \, F_{n \times n} \, D_{n \times n},
\end{equation}
where $S_{l \times n}$ is the matrix whose entries are all zeros,
aside from a single 1 in column $s_j$ of row $j$
for $j = 1$,~$2$, \dots, $l-1$,~$l$,
where $s_1$,~$s_2$, \dots, $s_{l-1}$,~$s_l$ are
i.i.d. integer random variables,
each distributed uniformly over $\{1, 2, \dots, n-1, n\}$;
moreover, $F_{n \times n}$ is the discrete Fourier transform
defined in~(\ref{DFT}),
and $D_{n \times n}$ is the diagonal matrix whose diagonal entries
$d_1$,~$d_2$, \dots, $d_{n-1}$,~$d_n$ are i.i.d. complex random variables,
each distributed uniformly over the unit circle.
(In our numerical implementations,
we drew $s_1$,~$s_2$, \dots, $s_{l-1}$,~$s_l$
from $\{1, 2, \dots, n-1, n\}$ without replacement,
instead of using i.i.d. draws.)
We observe that both $F_{n \times n}$ and $D_{n \times n}$ are unitary.

In~(\ref{SRFT}), $H_{n \times n}$ is the random matrix given by the formula
\begin{equation}
\label{Rokhlin_trans}
H_{n \times n} = \Theta_{n \times n} \, \Pi_{n \times n} \, Z_{n \times n}
              \, \tilde{\Theta}_{n \times n} \, \tilde{\Pi}_{n \times n}
              \, \tilde{Z}_{n \times n},
\end{equation}
where $\Pi_{n \times n}$ and $\tilde{\Pi}_{n \times n}$
are permutation matrices chosen independently and uniformly at random, and
$Z_{n \times n}$ and $\tilde{Z}_{n \times n}$ are diagonal matrices
whose diagonal entries $\zeta_1$,~$\zeta_2$, \dots, $\zeta_{n-1}$,~$\zeta_n$
and $\tilde{\zeta}_1$,~$\tilde{\zeta}_2$, \dots,
$\tilde{\zeta}_{n-1}$,~$\tilde{\zeta}_n$ are i.i.d. complex random variables,
each distributed uniformly over the unit circle; furthermore,
$\Theta_{n \times n}$ and $\tilde{\Theta}_{n \times n}$ are the matrices
defined via the formulae
\begin{multline}
\label{big_def}
\Theta_{n \times n} =
\left( \begin{array}{ccccc}
 \cos(\theta_1) & \sin(\theta_1) & 0 &      0 & 0 \\
-\sin(\theta_1) & \cos(\theta_1) & 0 &      0 & 0 \\
              0 &              0 & 1 &      0 & 0 \\
              0 &              0 & 0 & \ddots & 0 \\
              0 &              0 & 0 &      0 & 1
\end{array} \right) \cdot \\
\cdot \left( \begin{array}{ccccc}
 1 &               0 &              0 & 0 & 0 \\
 0 &  \cos(\theta_2) & \sin(\theta_2) & 0 & 0 \\
 0 & -\sin(\theta_2) & \cos(\theta_2) & 0 & 0 \\
 0 &               0 &              0 & 1 & 0 \\
 0 &               0 &              0 & 0 & \ddots
\end{array} \right) \cdots \\
\cdots \left( \begin{array}{ccccc}
\ddots & 0 &                   0 &                  0 & 0 \\
     0 & 1 &                   0 &                  0 & 0 \\
     0 & 0 &  \cos(\theta_{n-2}) & \sin(\theta_{n-2}) & 0 \\
     0 & 0 & -\sin(\theta_{n-2}) & \cos(\theta_{n-2}) & 0 \\
     0 & 0 &                   0 &                  0 & 1
\end{array} \right) \cdot \\
\cdot \left( \begin{array}{ccccc}
1 &      0 & 0 &                   0 &                  0 \\
0 & \ddots & 0 &                   0 &                  0 \\
0 &      0 & 1 &                   0 &                  0 \\
0 &      0 & 0 &  \cos(\theta_{n-1}) & \sin(\theta_{n-1}) \\
0 &      0 & 0 & -\sin(\theta_{n-1}) & \cos(\theta_{n-1})
\end{array} \right)
\end{multline}
and (the same as~(\ref{big_def}), but with tildes)
\begin{multline}
\tilde{\Theta}_{n \times n} =
\left( \begin{array}{ccccc}
 \cos(\tilde{\theta}_1) & \sin(\tilde{\theta}_1) & 0 &      0 & 0 \\
-\sin(\tilde{\theta}_1) & \cos(\tilde{\theta}_1) & 0 &      0 & 0 \\
                      0 &                      0 & 1 &      0 & 0 \\
                      0 &                      0 & 0 & \ddots & 0 \\
                      0 &                      0 & 0 &      0 & 1
\end{array} \right) \cdot \\
\cdot \left( \begin{array}{ccccc}
 1 &                       0 &                      0 & 0 & 0 \\
 0 &  \cos(\tilde{\theta}_2) & \sin(\tilde{\theta}_2) & 0 & 0 \\
 0 & -\sin(\tilde{\theta}_2) & \cos(\tilde{\theta}_2) & 0 & 0 \\
 0 &                       0 &                      0 & 1 & 0 \\
 0 &                       0 &                      0 & 0 & \ddots
\end{array} \right) \cdots \\
\cdots \left( \begin{array}{ccccc}
\ddots & 0 &                           0 &                          0 & 0 \\
     0 & 1 &                           0 &                          0 & 0 \\
     0 & 0 &  \cos(\tilde{\theta}_{n-2}) & \sin(\tilde{\theta}_{n-2}) & 0 \\
     0 & 0 & -\sin(\tilde{\theta}_{n-2}) & \cos(\tilde{\theta}_{n-2}) & 0 \\
     0 & 0 &                           0 &                          0 & 1
\end{array} \right) \cdot \\
\cdot \left( \begin{array}{ccccc}
1 &      0 & 0 &                           0 &                          0 \\
0 & \ddots & 0 &                           0 &                          0 \\
0 &      0 & 1 &                           0 &                          0 \\
0 &      0 & 0 &  \cos(\tilde{\theta}_{n-1}) & \sin(\tilde{\theta}_{n-1}) \\
0 &      0 & 0 & -\sin(\tilde{\theta}_{n-1}) & \cos(\tilde{\theta}_{n-1})
\end{array} \right),
\end{multline}
where $\theta_1$,~$\theta_2$, \dots, $\theta_{n-2}$,~$\theta_{n-1}$,
$\tilde{\theta}_1$,~$\tilde{\theta}_2$, \dots,
$\tilde{\theta}_{n-2}$,~$\tilde{\theta}_{n-1}$
are i.i.d. real random variables drawn uniformly from $[0,2 \pi]$.
We observe that $\Theta_{n \times n}$, $\tilde{\Theta}_{n \times n}$,
$\Pi_{n \times n}$, $\tilde{\Pi}_{n \times n}$,
$Z_{n \times n}$, and~$\tilde{Z}_{n \times n}$ are all unitary,
and so $H_{n \times n}$ is also unitary.

We call the transform $T_{l \times n}$ an ``SRFT'' for lack of a better term.

The following technical lemma is a slight reformulation
of Lemma~4.4 of~\cite{woolfe-liberty-rokhlin-tygert}.

\begin{lemma}
\label{reformulation}
Suppose that $\alpha$ and $\beta$ are real numbers greater than 1,
and $l$, $m$, and $n$ are positive integers, such that
\begin{equation}
\label{weak}
n > l \ge \frac{\alpha^2 \, \beta}{(\alpha-1)^2} \, m^2.
\end{equation}
Suppose further that $T_{l \times n}$ is the SRFT defined in~(\ref{SRFT}),
and that $Q_{n \times m}$ is a matrix whose columns are orthonormal.

Then, the least singular value $\sigma_m$ of $T_{l \times n} \, Q_{n \times m}$
satisfies
\begin{equation}
\label{inverse}
\sigma_m \ge \sqrt{\frac{l}{\alpha n}}
\end{equation}
with probability at least $1-\frac{1}{\beta}$.
\end{lemma}

\section{Description of the algorithm}
\label{algorithm}

Suppose that $\epsilon$ is a positive real number less than 1,
and $l$, $m$, and $n$ are positive integers with $m < l < n$.
Suppose further that $A_{m \times n}$ is a full-rank matrix,
$b_{m \times 1}$ is a vector, and $p_{n \times 1}$ is the vector
of minimal Euclidean norm satisfying
$A_{m \times n} \, p_{n \times 1} = b_{m \times 1}$.
In order to construct a vector $x_{n \times 1}$ such that
$\| x_{n \times 1} - p_{n \times 1} \| \le \epsilon \, \| p_{n \times 1} \|$
with high probability (increasingly high probability as a parameter $\alpha>1$
increases), we compute the vector $c_{n \times 1}$ of minimal Euclidean norm
that is a linear combination of random vectors
and satisfies $A_{m \times n} \, c_{n \times 1} = b_{m \times 1}$,
then use the algorithm of~\cite{rokhlin-tygert}
to compute the orthogonal projection of $c_{n \times 1}$
onto the column span of $(A_{m \times n})^*$.
More precisely, we perform the following five steps:

\begin{enumerate}
\item Construct the matrix
      \begin{equation}
      \label{randmat}
      S_{l \times m} = T_{l \times n} \, (A_{m \times n})^*,
      \end{equation}
      applying the SRFT $T_{l \times n}$ defined in~(\ref{SRFT})
      to every column of $(A_{m \times n})^*$
      (see, for example, Subsection~3.3 of~\cite{woolfe-liberty-rokhlin-tygert}
      for details on applying the SRFT rapidly).
\item Construct the vector $z_{l \times 1}$ of minimal Euclidean norm
      solving the system of linear equations
      \begin{equation}
      \label{lineq}
      (S_{l \times m})^* \, z_{l \times 1} = b_{m \times 1},
      \end{equation}
      where $S_{l \times m}$ is the matrix defined in~(\ref{randmat})
      (see, for example, Algorithm~5.7.2 in~\cite{golub-van_loan}
      for details on constructing $z_{l \times 1}$).
\item Construct the vector
      \begin{equation}
      \label{nonminnorm}
      c_{n \times 1} = (T_{l \times n})^* \, z_{l \times 1},
      \end{equation}
      where $z_{l \times 1}$ is the vector of minimal Euclidean norm
      solving~(\ref{lineq}), and $T_{l \times n}$ is the same realization
      of the SRFT as in~(\ref{randmat})
      (see, for example, Subsection~3.3 of~\cite{woolfe-liberty-rokhlin-tygert}
      for details on applying the adjoint of the SRFT rapidly).
\item Use the algorithm of~\cite{rokhlin-tygert} for the construction
      of a vector $y_{m \times 1}$ minimizing
      \begin{equation}
      \label{overdet}
      \| (A_{m \times n})^* \, y_{m \times 1} - c_{n \times 1} \|^2
      \end{equation}
      to relative precision $\epsilon^2 l/(\alpha n)$ or better,
      where $c_{n \times 1}$ is the vector defined in~(\ref{nonminnorm}).
      The parameter $l$ for the algorithm of~\cite{rokhlin-tygert}
      should be the same as for the present algorithm.
\item Construct the desired vector
      \begin{equation}
      \label{output}
      x_{n \times 1} = (A_{m \times n})^* \, y_{m \times 1},
      \end{equation}
      where $y_{m \times 1}$ is the vector from Step~4.
\end{enumerate}

\begin{remark}
In Step~2 above, we assume that $S_{l \times m}$ defined in~(\ref{randmat})
is a full-rank matrix. Lemma~\ref{reformulation} in Section~\ref{prelims}
above guarantees this with high probability when $l > m^2$,
by taking $\alpha$ in the lemma arbitrarily large;
numerical experiments indicate that $l > m$ suffices.
\end{remark}

\begin{remark}
It is possible to improve the approximation $x_{n \times 1}$
via preconditioned conjugate gradient iterations similar to those
proposed in~\cite{rokhlin-tygert}.
However, the approximation produced by the above algorithm
is already highly accurate (see, for example, Section~\ref{numerical}
or Theorem~\ref{main_theorem} in Section~\ref{proofs} below),
and further iterative improvement may double the running time of the algorithm.
\end{remark}

\section{Proof of accuracy}
\label{proofs}

In this section, we prove Theorem~\ref{main_theorem},
guaranteeing that the algorithm of Section~\ref{algorithm}
produces high accuracy with high probability.

The following lemma states that the orthogonal projection
onto the column span of $(A_{m \times n})^*$
of the vector $c_{n \times 1}$ defined in~(\ref{nonminnorm})
is the vector of minimal Euclidean norm that the algorithm aims to approximate.

\begin{lemma}
\label{baselemma}
Suppose that $l$, $m$, and $n$ are positive integers with $m < l < n$.
Suppose further that $A_{m \times n}$ is a full-rank matrix,
$b_{m \times 1}$ is a vector,
$S_{l \times m}$ is the matrix defined in~(\ref{randmat})
and is a full-rank matrix,
and $c_{n \times 1}$ is the vector defined in~(\ref{nonminnorm}).

Then, the orthogonal projection $p_{n \times 1}$
of $c_{n \times 1}$ onto the column span of $(A_{m \times n})^*$
is the vector of minimal Euclidean norm satisfying
\begin{equation}
\label{basicproj}
A_{m \times n} \, p_{n \times 1} = b_{m \times 1}.
\end{equation}
\end{lemma}

\begin{proof}
Combining~(\ref{randmat}), (\ref{lineq}), and~(\ref{nonminnorm}) yields that
\begin{equation}
\label{easy}
A_{m \times n} \, c_{n \times 1} = b_{m \times 1}.
\end{equation}
Combining~(\ref{easy}) and the fact that $p_{n \times 1}$
is the orthogonal projection of $c_{n \times 1}$ onto the column span
of $(A_{m \times n})^*$ completes the proof.
\end{proof}

The following lemma states that, with high probability, the Euclidean norm
of the vector $c_{n \times 1}$ defined in~(\ref{nonminnorm})
is not too much greater than the Euclidean norm
of its orthogonal projection $p_{n \times 1}$ onto the column span
of $(A_{m \times n})^*$.

\begin{lemma}
\label{problemma}
Suppose that $\alpha$ and $\beta$ are real numbers greater than 1,
and $l$, $m$, and $n$ are positive integers, such that~(\ref{weak}) holds.
Suppose further that $A_{m \times n}$ is a full-rank matrix,
$b_{m \times 1}$ is a vector,
$c_{n \times 1}$ is the vector defined in~(\ref{nonminnorm}),
and $p_{n \times 1}$ is the orthogonal projection of $c_{n \times 1}$
onto the column span of $(A_{m \times n})^*$.

Then,
\begin{equation}
\label{energy}
\| c_{n \times 1} \| \le \sqrt{\frac{\alpha n}{l}} \, \| p_{n \times 1} \|
\end{equation}
with probability at least $1-\frac{1}{\beta}$.
\end{lemma}

\begin{proof}
Using the fact that $A_{m \times n}$ is a full-rank matrix,
we construct a $QR$-decomposition
\begin{equation}
\label{QR}
(A_{m \times n})^* = Q_{n \times m} \, R_{m \times m}
\end{equation}
such that the columns of $Q_{n \times m}$ are an orthonormal basis
for the column span of $(A_{m \times n})^*$.
We first show that the SRFT $T_{l \times n}$ used
in~(\ref{randmat}) and~(\ref{nonminnorm}) provides
\begin{equation}
\label{intermed}
\| z_{l \times 1} \|
\le \sqrt{\frac{\alpha n}{l}}
 \, \| (Q_{n \times m})^* \, (T_{l \times n})^* \, z_{l \times 1} \|
\end{equation}
with probability at least $1-\frac{1}{\beta}$,
where $z_{l \times 1}$ is the vector of minimal Euclidean norm
solving~(\ref{lineq}).
We then express the left- and right-hand sides of~(\ref{intermed})
in terms of $c_{n \times 1}$ and $p_{n \times 1}$,
rather than $z_{l \times 1}$, in order to obtain~(\ref{energy}).

It follows from the fact that $z_{l \times 1}$ is the vector
of minimal Euclidean norm solving~(\ref{lineq})
that $z_{l \times 1}$ belongs to the column span
of $S_{l \times m}$ from~(\ref{lineq}), that is,
there exists a vector $w_{m \times 1}$ such that
\begin{equation}
\label{linearcombo}
z_{l \times 1} = S_{l \times m} \, w_{m \times 1}.
\end{equation}
Combining~(\ref{linearcombo}), (\ref{randmat}), and~(\ref{QR}) yields that
\begin{equation}
\label{uno}
\| z_{l \times 1} \|^2 = (z_{l \times 1})^* \, z_{l \times 1}
= (w_{m \times 1})^* \, (R_{m \times m})^* \, (Q_{n \times m})^* \,
  (T_{l \times n})^* \, z_{l \times 1}.
\end{equation}
The Cauchy-Schwarz inequality yields that
\begin{equation}
\label{dos}
(w_{m \times 1})^* \, (R_{m \times m})^* \, (Q_{n \times m})^* \,
(T_{l \times n})^* \, z_{l \times 1}
\le \| R_{m \times m} \, w_{m \times 1} \|
\, \| (Q_{n \times m})^* \, (T_{l \times n})^* \, z_{l \times 1} \|.
\end{equation}
It follows from~(\ref{inverse}) that
\begin{equation}
\label{tres}
\| R_{m \times m} \, w_{m \times 1} \|
\le \sqrt{\frac{\alpha n}{l}} \,
\| T_{l \times n} \, Q_{n \times m} \, R_{m \times m} \, w_{m \times 1} \|
\end{equation}
with probability at least $1-\frac{1}{\beta}$.
Combining~(\ref{QR}), (\ref{randmat}), and~(\ref{linearcombo}) yields that
\begin{equation}
\label{cuatro}
T_{l \times n} \, Q_{n \times m} \, R_{m \times m} \, w_{m \times 1}
= z_{l \times 1}.
\end{equation}
Combining~(\ref{uno}), (\ref{dos}), (\ref{tres}), and~(\ref{cuatro})
yields~(\ref{intermed}).

We now express the left- and right-hand sides of~(\ref{intermed})
in terms of $c_{n \times 1}$ and $p_{n \times 1}$,
rather than $z_{l \times 1}$.

First, we consider the left-hand side of~(\ref{intermed}).
Combining~(\ref{nonminnorm}) and the fact that the columns
of $(T_{l \times n})^*$ are orthonormal yields that
\begin{equation}
\label{intermed1}
\| z_{l \times 1} \| = \| c_{n \times 1} \|.
\end{equation}

Next, we consider the right-hand side of~(\ref{intermed}).
It follows from the fact that the columns of $Q_{n \times m}$
are an orthonormal basis for the column span of $(A_{m \times n})^*$
that the orthogonal projection $p_{n \times 1}$
of $c_{n \times 1}$ onto the column span of $(A_{m \times n})^*$ is
\begin{equation}
\label{projector}
p_{n \times 1} = Q_{n \times m} \, (Q_{n \times m})^* \, c_{n \times 1}.
\end{equation}
Combining~(\ref{projector}) and the fact that the columns
of $Q_{n \times m}$ are orthonormal yields that
\begin{equation}
\label{normproj}
\| p_{n \times 1} \| = \| (Q_{n \times m})^* \, c_{n \times 1} \|.
\end{equation}
Combining~(\ref{normproj}) and~(\ref{nonminnorm}) yields that
\begin{equation}
\label{intermed2}
\| p_{n \times 1} \|
= \| (Q_{n \times m})^* \, (T_{l \times n})^* \, z_{l \times 1}\|.
\end{equation}

Finally, combining~(\ref{intermed}), (\ref{intermed1}), and~(\ref{intermed2})
yields~(\ref{energy}).
\end{proof}

The following lemma states that the vector $x_{n \times 1}$
produced by the algorithm is an accurate approximation
to the orthogonal projection onto the column span of $(A_{m \times n})^*$
of the vector $c_{n \times 1}$ defined in~(\ref{nonminnorm}),
provided that the projection $p_{n \times 1}$ satisfies~(\ref{energy}).

\begin{lemma}
\label{mainlemma}
Suppose that $\epsilon$ and $\alpha$ are positive real numbers
with $\epsilon < 1 < \alpha$,
and $l$, $m$, and $n$ are positive integers with $m < l < n$.
Suppose further that $A_{m \times n}$ is a full-rank matrix,
$b_{m \times 1}$ is a vector,
$S_{l \times m}$ is the matrix defined in~(\ref{randmat})
and is a full-rank matrix,
$c_{n \times 1}$ is the vector defined in~(\ref{nonminnorm}),
$p_{n \times 1}$ is the orthogonal projection of $c_{n \times 1}$
onto the column span of $(A_{m \times n})^*$,
$y_{m \times 1}$ is a vector minimizing~(\ref{overdet})
to relative precision $\epsilon^2 l/(\alpha n)$ or better,
and $x_{n \times 1}$ is the vector defined in~(\ref{output}).
Suppose in addition that~(\ref{energy}) holds.

Then,
\begin{equation}
\label{goodapprox}
\| x_{n \times 1} - p_{n \times 1} \|
\le \epsilon \, \| p_{n \times 1} \|.
\end{equation}
\end{lemma}

\begin{proof}
It follows from~(\ref{output}) that $x_{n \times 1}$
belongs to the column span of $(A_{m \times n})^*$.
Combining this fact and the fact that $c_{n \times 1} - p_{n \times 1}$
is the orthogonal projection of $c_{n \times 1}$
onto the orthogonal complement of the column span of $(A_{m \times n})^*$
yields that $c_{n \times 1} - p_{n \times 1}$ is the orthogonal projection
of $c_{n \times 1} - x_{n \times 1}$
onto the orthogonal complement of the column span of $(A_{m \times n})^*$.
Similarly, $p_{n \times 1} - x_{n \times 1}$ is the orthogonal projection
of $c_{n \times 1} - x_{n \times 1}$
onto the column span of $(A_{m \times n})^*$.
We thus obtain the Pythagorean identity
\begin{equation}
\label{first}
\| c_{n \times 1} - p_{n \times 1} \|^2
+ \| p_{n \times 1} - x_{n \times 1} \|^2
= \| c_{n \times 1} - x_{n \times 1} \|^2.
\end{equation}

It follows from the fact that $p_{n \times 1}$ is the orthogonal projection
of $c_{n \times 1}$ onto the column span of $(A_{m \times n})^*$ that
the minimal value of~(\ref{overdet})
is $\| p_{n \times 1} - c_{n \times 1} \|^2$.
Combining this fact, (\ref{output}), and the fact that $y_{m \times 1}$
minimizes~(\ref{overdet}) to relative precision $\epsilon^2 l/(\alpha n)$
or better yields that
\begin{equation}
\label{second}
\| x_{n \times 1} - c_{n \times 1} \|^2-\| p_{n \times 1} - c_{n \times 1} \|^2
\le \frac{\epsilon^2 l}{\alpha n} \, \| p_{n \times 1} - c_{n \times 1} \|^2.
\end{equation}

Combining~(\ref{first}) and~(\ref{second}) yields that
\begin{equation}
\label{third}
\| x_{n \times 1} - p_{n \times 1} \|^2
\le \frac{\epsilon^2 l}{\alpha n} \, \| c_{n \times 1} - p_{n \times 1} \|^2.
\end{equation}

It follows from the fact that $c_{n \times 1} - p_{n \times 1}$
is the orthogonal projection of $c_{n \times 1}$
onto the orthogonal complement of the column span of $(A_{m \times n})^*$ that
\begin{equation}
\label{fourth}
\| c_{n \times 1} - p_{n \times 1} \| \le \| c_{n \times 1} \|.
\end{equation}
Combining~(\ref{third}), (\ref{fourth}), and~(\ref{energy})
yields~(\ref{goodapprox}).
\end{proof}

Combining Lemmas~\ref{reformulation}, \ref{baselemma}, \ref{problemma},
and \ref{mainlemma} yields the following theorem,
guaranteeing that the algorithm produces high accuracy with high probability.

\begin{theorem}
\label{main_theorem}
Suppose that $\epsilon$, $\alpha$, and $\beta$ are positive real numbers
with $\epsilon < 1 < \alpha$ and $\beta > 2$,
and $l$, $m$, and $n$ are positive integers, such that~(\ref{weak}) holds.
Suppose further that $A_{m \times n}$ is a full-rank matrix,
$b_{m \times 1}$ is a vector,
and $p_{n \times 1}$ is the vector of minimal Euclidean norm satisfying
\begin{equation}
A_{m \times n} \, p_{n \times 1} = b_{m \times 1}.
\end{equation}

Then, the vector $x_{n \times 1}$ defined in~(\ref{output}) satisfies
\begin{equation}
\| x_{n \times 1} - p_{n \times 1} \| \le \epsilon \, \| p_{n \times 1} \|
\end{equation}
with probability at least $1-\frac{2}{\beta}$.
\end{theorem}

\begin{remark}
The probability of failure in Theorem~\ref{main_theorem} is $2/\beta$,
rather than just the probability $1/\beta$ of failure
in Lemma~\ref{problemma}, in order to account for the possibility that
the algorithm of~\cite{rokhlin-tygert} fails to produce a vector
$y_{m \times 1}$ minimizing~(\ref{overdet}) to relative precision
$\epsilon^2 l/(\alpha n)$ or better,
after using at most the number of floating-point operations
specified in Section~\ref{costs} below.
\end{remark}

\begin{remark}
Empirically, requiring~(\ref{weak}) appears to be excessive.
In practice, choosing any $l \ge 4m$ and $\alpha = 4$ makes failure
of the algorithm too improbable to detect
({\it cf.} Remark~2 in~\cite{rokhlin-tygert}).
\end{remark}

\section{Computational costs}
\label{costs}

In this section, we tabulate the numbers of floating-point operations
required by the five steps in the algorithm of Section~\ref{algorithm}:

\begin{enumerate}
\item Applying $T_{l \times n}$ to every column of $(A_{m \times n})^*$
      costs $\bigoh(mn \, \log(l))$.
\item Constructing $z_{l \times 1}$ costs $\bigoh(l m^2)$.
\item Applying $(T_{l \times n})^*$ to $z_{l \times 1}$
      costs $\bigoh(n \, \log(n))$.
\item Constructing $y_{l \times 1}$ via the algorithm of~\cite{rokhlin-tygert}
      costs $\bigoh(mn \, \log(\alpha n/\epsilon^2) + l m^2)$.
\item Applying $(A_{m \times n})^*$ to $y_{m \times 1}$ costs $\bigoh(mn)$.
\end{enumerate}

Summing up the costs in the five steps above,
and using the facts that $m < l < n$ and $\epsilon < 1 < \alpha$,
we see that the cost of the entire algorithm is
\begin{equation}
\label{theoretical}
C_{\rm theoretical} = \bigoh(mn \, \log(\sqrt{\alpha n}/\epsilon) + l m^2)
\end{equation}
floating-point operations.
Theorem~\ref{main_theorem} in Section~\ref{proofs} above
guarantees that the algorithm produces
high accuracy with high probability when $l$ and $\alpha$ satisfy~(\ref{weak}).
In practice, choosing $l = 4m$ and $\alpha = 4$ makes failure
of the algorithm too improbable to detect
(see also Remark~2 in~\cite{rokhlin-tygert}),
and the cost in~(\ref{theoretical}) then becomes
\begin{equation}
C_{\rm typical} = \bigoh(mn \, \log(\sqrt{n}/\epsilon) + m^3)
\end{equation}
floating-point operations.

\section{Numerical Results}
\label{numerical}

In this section, we describe the results of several numerical tests
of the algorithm of the present paper.

For various positive integers $m$ and $n$ with $m < n$,
we use the algorithm to compute a vector $x_{n \times 1}$
approximating the vector $p_{n \times 1}$ of minimal Euclidean norm
satisfying $A_{m \times n} \, p_{n \times 1} = b_{m \times 1}$,
where $b_{m \times 1}$ is a vector,
and $A_{m \times n}$ is the matrix defined via the formula
\begin{equation}
\label{matrix_definition}
A_{m \times n} = U_{m \times m} \, \Sigma_{m \times m} \, (V_{n \times m})^*;
\end{equation}
in the experiments described below,
$U_{m \times m}$ is obtained by applying the Gram-Schmidt process
to the columns of an $m \times m$ matrix whose entries are
i.i.d. centered complex Gaussian random variables,
$V_{n \times m}$ is obtained by applying the Gram-Schmidt process
to the columns of an $n \times m$ matrix whose entries are
i.i.d. centered complex Gaussian random variables,
and $\Sigma_{m \times m}$ is a diagonal matrix, with the diagonal entries
\begin{equation}
\Sigma_{j,j} = 10^{-6(j-1)/(m-1)}
\end{equation}
for $j = 1$,~$2$, \dots, $m-1$,~$m$.
Clearly, the condition number $\kappa_A$ of $A_{m \times n}$ is
\begin{equation}
\label{condnum}
\kappa_A = \Sigma_{1,1}/\Sigma_{m,m} = 10^6.
\end{equation}
The vector $b_{m \times 1}$ is defined via the formula
\begin{equation}
\label{rhs_definition}
b_{m \times 1} = A_{m \times n} \, p_{n \times 1},
\end{equation}
where $p_{n \times 1}$ is the vector defined via the formula
\begin{equation}
\label{exact_solution}
p_{n \times 1}
= \frac{1}{\sqrt{m}} \sum_{j=1}^m \epsilon_j \, v_{n \times 1}^{(j)},
\end{equation}
with $\epsilon_1$,~$\epsilon_2$, \dots, $\epsilon_{m-1}$,~$\epsilon_m$
being pseudorandom positive and negative ones,
and $v_{n \times 1}^{(1)}$,~$v_{n \times 1}^{(2)}$, \dots,
$v_{n \times 1}^{(m-1)}$,~$v_{n \times 1}^{(m)}$ being the columns
of $V_{n \times m}$.

\begin{remark}
By construction, $p_{n \times 1}$ is the vector of minimal Euclidean norm
such that $A_{m \times n} \, p_{n \times 1} = b_{m \times 1}$;
the Euclidean norm of $p_{n \times 1}$ is minimal since $p_{n \times 1}$
belongs to the column span of $(A_{m \times n})^*$.
The aim of the algorithm is to construct an approximation $x_{n \times 1}$
to $p_{n \times 1}$.
\end{remark}

For the direct computations, we used the classical algorithm
for pivoted $QR$-de\-com\-po\-si\-tions
based on plane (Householder) reflections
(see, for example, Chapter~5 in~\cite{golub-van_loan}).
We implemented the algorithms in Fortran~77 in double-precision arithmetic,
and used the Lahey/Fujitsu Express v6.2 compiler,
with the optimization flag {\tt -{}-o2} enabled.
We used one core of a 1.86~GHz Intel Centrino Core Duo microprocessor
with 2~MB of L2 cache and 1~GB of RAM.
For the algorithm of~\cite{rokhlin-tygert} used in Step~4
of the algorithm of Section~\ref{algorithm},
we requested that $y_{m \times 1}$ minimize~(\ref{overdet})
to relative precision $(10^{-14} \cdot \kappa_A)^2 \cdot m/n$ or better,
where $\kappa_A$ is the condition number of $A_{m \times n}$
given in~(\ref{condnum}).
We used a double-precision version of P.~N. Swarztrauber's FFTPACK library
for fast Fourier transforms.

Table~1 displays timing results with $n = 16384$ for various values of $m$;
Table~2 displays the corresponding errors.
Table~3 displays timing results with $m = 256$ for various values of $n$;
Table~4 displays the corresponding errors.

\begin{table}
\hfil {\bf Table 1} \hfil\hfil {\bf Table 2} \hfil \\\vspace{-.75em}\\
\begin{tabular}{c|c|c|c|c|c}
$m$ &   $n$ &  $l$ & $t_0$ & $t_{\rm r}$ & $t_0/t_{\rm r}$ \\\hline
128 & 16384 &  512 & .27E1 &       .24E1 &             1.2 \\
256 & 16384 & 1024 & .11E2 &       .56E1 &             2.0 \\
512 & 16384 & 2048 & .60E2 &       .20E2 &             3.0 \\
\end{tabular}
\hfil
\begin{tabular}{c|c|c|c|c}
$m$ &   $n$ &  $l$ & $\epsilon_0$ & $\epsilon_{\rm r}$ \\\hline
128 & 16384 &  512 &     .14E--14 &           .16E--14 \\
256 & 16384 & 1024 &     .11E--14 &           .17E--14 \\
512 & 16384 & 2048 &     .80E--15 &           .29E--14 \\
\end{tabular}

\bigskip
\bigskip

\hfil {\bf Table 3} \hfil\hfil {\bf Table 4} \hfil \\\vspace{-.75em}\\
\begin{tabular}{c|c|c|c|c|c}
$m$ &   $n$ &  $l$ & $t_0$ & $t_{\rm r}$ & $t_0/t_{\rm r}$ \\\hline
256 &  4096 & 1024 & .26E1 &       .20E1 &             1.3 \\
256 &  8192 & 1024 & .51E1 &       .32E1 &             1.6 \\
256 & 16384 & 1024 & .11E2 &       .56E1 &             2.0 \\
256 & 32768 & 1024 & .29E2 &       .16E2 &             2.5 \\
\end{tabular}
\hfil
\begin{tabular}{c|c|c|c|c}
$m$ &   $n$ &  $l$ & $\epsilon_0$ & $\epsilon_{\rm r}$ \\\hline
256 &  4096 & 1024 &     .27E--15 &           .31E--14 \\
256 &  8192 & 1024 &     .45E--15 &           .27E--14 \\
256 & 16384 & 1024 &     .11E--14 &           .17E--14 \\
256 & 32768 & 1024 &     .22E--14 &           .16E--14 \\
\end{tabular}

\bigskip
\end{table}

The headings of the tables have the following meanings:
\begin{itemize}
\item $m$ is the number of rows in the matrix $A_{m \times n}$,
      as well as the length of the vector $b_{m \times 1}$,
      in $A_{m \times n} \, p_{n \times 1} = b_{m \times 1}$.
\item $n$ is the number of columns in the matrix $A_{m \times n}$,
      as well as the length of the vector $p_{n \times 1}$,
      in $A_{m \times n} \, p_{n \times 1} = b_{m \times 1}$.
\item $l$ is the number of rows in the matrix $T_{l \times n}$ used
      in Steps~1 and~3 of the algorithm of Section~\ref{algorithm},
      as well as the analogous parameter used in the algorithm
      of~\cite{rokhlin-tygert} needed in Step~4.
\item $t_0$ is the time in seconds required
      by the direct, classical algorithm.
\item $t_{\rm r}$ is the time in seconds required
      by the algorithm of the present paper.
\item $t_0/t_{\rm r}$ is the factor by which
      the algorithm of the present paper
      is faster than the classical algorithm.
\item $\epsilon_0$ is defined via the formula
\begin{equation}
\epsilon_0 = \frac{\| x_{n \times 1}^{(0)} - p_{n \times 1} \|}
                  {\kappa_A \, \| p_{n \times 1} \|},
\end{equation}
where $\kappa_A$ is the condition number of $A_{m \times n}$ given
in~(\ref{condnum}), and $x_{n \times 1}^{(0)}$ is the vector
produced by the direct, classical algorithm
approximating the vector $p_{n \times 1}$ of minimal Euclidean norm
such that $A_{m \times n} \, p_{n \times 1} = b_{m \times 1}$.
\item $\epsilon_{\rm r}$ is defined via the formula
\begin{equation}
\epsilon_{\rm r} = \frac{\| x_{n \times 1} - p_{n \times 1} \|}
                        {\kappa_A \, \| p_{n \times 1} \|},
\end{equation}
where $\kappa_A$ is the condition number of $A_{m \times n}$ given
in~(\ref{condnum}), and $x_{n \times 1}$ is the vector
produced by the algorithm of Section~\ref{algorithm}
approximating the vector $p_{n \times 1}$ of minimal Euclidean norm
such that $A_{m \times n} \, p_{n \times 1} = b_{m \times 1}$.
\end{itemize}

\begin{remark}
Standard perturbation theory shows that
$\epsilon_0$ and $\epsilon_{\rm r}$ are the appropriately normalized measures
of the relative precision produced by the algorithms;
see, for example, Section~5.5.3 in~\cite{dahlquist-bjorck}.
\end{remark}

The values for $\epsilon_{\rm r}$ reported in the tables are
the worst (maximal) values encountered during 10 independent randomized trials
of the algorithm, as applied to the same matrix $A_{m \times n}$
and vector $b_{m \times 1}$.
The values for $t_{\rm r}$ reported in the tables are the average values
over 10 independent randomized trials.
None of the quantities reported in the tables varied significantly
over repeated randomized trials.

The following observations can be made from the examples reported here,
and from our more extensive experiments:

\begin{enumerate}
\item When $m = 512$, $n = 16384$,
and the condition number of $A_{m \times n}$ is $10^6$,
the randomized algorithm runs 3 times faster than the classical algorithm
based on plane (Householder) reflections, even at full double precision.
\item Our choice $l = 4m$ appears to make failure of the algorithm
too improbable to detect.
\item The algorithm of the present article reliably produces high precision
at reasonably low cost.
\end{enumerate}

\section{Conclusion}
\label{conclusions}

This article provides a fast algorithm for computing the minimal-norm solution
to an underdetermined system of linear equations.
If the matrices $A_{m \times n}$ and $(A_{m \times n})^*$
associated with the system of linear equations
can be applied sufficiently rapidly to arbitrary vectors, then the algorithm
of the present paper can be accelerated further.

The theoretical bounds in Lemma~\ref{reformulation}, Lemma~\ref{problemma},
and Theorem~\ref{main_theorem} should be considered preliminary.
Our numerical experiments indicate that the algorithm
of the present article performs better than our estimates guarantee.
Furthermore, there is nothing magical
about the subsampled randomized Fourier transform defined in~(\ref{SRFT}).
In our experience, several other similar transforms seem
to work at least as well, and we are investigating these alternatives
(see, for example,~\cite{ailon-liberty}).

\section*{Acknowledgments}
We would like to thank Vladimir Rokhlin and Arthur Szlam
for helpful discussions.

\newpage

\bibliographystyle{siam}
\bibliography{under}

\end{document}